\documentclass{amsart}
\usepackage{amsmath,amssymb, amsbsy}
\usepackage[dvips]{graphicx}
\newcommand{\R}{{\mathbb R}}
\newcommand{\N}{{\mathbb N}}
\newcommand{\Z}{{\mathbb Z}}
\newcommand{\C}{{\mathbb C}}
\newcommand{\SN}{{\mathbb S}^{N-1}}

\newcommand{\e }{\varepsilon}
\newcommand{\Di}{{\mathcal D}^{1,2}(\R^N,\C)}

\newcommand{\A}{{\mathbf A}}

\renewcommand{\ge }{\geqslant}
\renewcommand{\geq }{\geqslant}

\renewcommand{\leq }{\leqslant}

\newenvironment{pf}{\noindent{\sc Proof}.\enspace}{\hfill\qed\medskip}
\newenvironment{pfn}[1]{\noindent{\bf Proof of
    {#1}.\enspace}}{\hfill\qed\medskip}
\newtheorem{Theorem}{Theorem}[section]
\newtheorem{Corollary}[Theorem]{Corollary}
\newtheorem{Lemma}[Theorem]{Lemma}
\newtheorem{Proposition}[Theorem]{Proposition}
\theoremstyle{definition}

\begin{document}

\title[Schr\"odinger equations with singular electromagnetic
potential]{A note on local asymptotics of solutions to singular elliptic equations via 
monotonicity methods}

\author[Veronica Felli \and Alberto Ferrero \and Susanna
Terracini]{Veronica Felli \and Alberto Ferrero \and Susanna Terracini}
\address{\hbox{\parbox{5.7in}{\medskip\noindent{Universit\`a di Milano
        Bicocca,\\
        Dipartimento di Ma\-t\-ema\-ti\-ca e Applicazioni, \\
        Via Cozzi
        53, 20125 Milano, Italy. \\[3pt]
        \em{E-mail addresses: }{\tt veronica.felli@unimib.it,
          alberto.ferrero@unimib.it, susanna.terracini@unimib.it}.}}}}

\date{}


 \begin{abstract}
   \noindent This paper completes and partially improves some of the
   results of \cite{FFT} about the asymptotic behavior of solutions of
   linear and nonlinear elliptic equations with singular coefficients
   via an Almgren type monotonicity formula
\end{abstract}

\maketitle

\section{Introduction and main results}
Regularity properties of solutions to linear elliptic partial
differential equations have been widely studied in the literature,
both in the case of singular coefficients in the elliptic operator and
in the case of domains with non smooth boundary. In order to determine
the regularity of solutions, some authors found proper asymptotic
expansions near the singularity of the coefficients or near a non
regular point of the boundary, see
\cite{CoDa,CoDaDu,FFT,MaMc,MaNaPl-1,MaNaPl-2,Mazzeo-1,Mazzeo-2} and the
references therein.

Our paper \cite{FFT} is concerned with the asymptotic
behavior near the singularity of solutions to equations associated
to the following class of Schr\"odinger operators with singular
homogeneous electromagnetic potentials:
\begin{equation*}
  {\mathcal L}_{\A,a}:=
  \left(-i\,\nabla+\frac{\A\big(\frac{x}{|x|}\big)}
    {|x|}\right)^{\!\!2}-\dfrac{a\big(\frac{x}{|x|}\big)}{|x|^2}.
\end{equation*}
In \cite{FFT}, we study both linear and nonlinear equations obtained as
perturbations of the operator ${\mathcal L}_{\A,a}$ in a domain
$\Omega\subset \R^N$ containing either the origin or a
neighborhood of $\infty$. More precisely, we deal with linear
equations of the type
\begin{equation} \label{u}
\mathcal L_{\A,a} u=h(x)\, u, \qquad {\rm in} \ \Omega,
\end{equation}
where $h\in L^\infty_{{\rm loc}}(\Omega\setminus \{0\})$ is
negligible with respect to the inverse square potential $|x|^{-2}$
near the singularity, and semilinear equations
\begin{equation} \label{nonlin0}
\mathcal L_{\A,a}u(x)=f(x,u(x))
\end{equation}
with $f$ having at most critical growth.  By solutions of \eqref{u} or
\eqref{nonlin0} we mean functions which belong to a suitable Sobolev
space depending on the magnetic potential $\A$ and solve the
corresponding equations in a distributional sense.

As far as the linear equation \eqref{u} is concerned, the main result
of \cite{FFT} provides the leading term in the asymptotic expansion
near the singularity of the coefficients. Similar
 asymptotic expansions were proved by Mazzeo \cite{Mazzeo-1},
\cite{Mazzeo-2}, with a completely different approach, in the more
general setting of elliptic equations on compact manifolds with
boundary (see also \cite{lesch}, \cite{LoMc}, and \cite{pacard}).

The main novelty of our approach in \cite{FFT} is the use of the
Almgren's monotonicity formula \cite{almgren}. This was the approach
already adopted in earlier works by Garofalo and Lin and then followed
by Kurata in \cite{Kurata} in order to prove the unique continuation
property.

In the present paper we illustrate the strengths of the monotonicity
formula approach, by completing and improving some of the results
obtained in \cite{FFT}.  The main purposes of this note are
essentially the following:
\begin{itemize}
\item[-] to deduce from the monotonicity formula more precise
  informations on the first term in the asymptotic expansion of
  \cite{Mazzeo-1}, \cite{Mazzeo-2} under some alternative assumptions
  on the perturbation~$h$ which require some integrability type conditions
  instead of pointwise decay as in~\cite{FFT},
\item[-] to provide a general method with the perspective of unifying the
  approach to linear and nonlinear equations with singular
  coefficients,
\item[-] to improve in the nonlinear case the results that in \cite{FFT}
  were obtained by using a-priori pointwise estimates on solutions.
\end{itemize}
In the remaining part of the introduction we will examine these three
goals with more detail.

Let us introduce some notations taken from
\cite{Mazzeo-1} adapting them to our context. Let us consider the
case where $\Omega=B_R=\{x\in\R^N:|x|<R\}$ for some $R>0$ in such a way that
$\overline\Omega\setminus \{0\}$ may be identified with the
cylinder $\SN\times (0,R]\subset \R^{N+1}$. If we identify the set
$\SN\times \{0\}$ to a point through an equivalence relation
$\sim$, then the quotient topological space $X:=(\SN\times
[0,R])/\sim$ becomes homeomorphic to $\overline B_R$. The topological space
$X$ has a natural structure of a compact manifold with boundary
$\partial X$ homeomorphic to $\SN$. On $X$ we can use the polar
coordinates $(r,\theta)$ with $r\in[0,R]$ and $\theta\in \SN$. If
we introduce the metric $g=dr^2+r^2 g_{\SN}$, where $g_{\SN}$ is
the standard metric on the unit sphere, then $X$ becomes a
Riemannian manifold isometric to $\overline B_R$. According to the
definition and the notations of \cite{Mazzeo-1}, a second order
elliptic operator $\mathcal L$ on $X$ is an operator which admits
a representation with respect to the coordinates $(r,\theta)$ of
the type
\begin{equation} \label{elliptic-Mazzeo}
\mathcal L=\sum_{0\leq j+|\beta|\leq 2} a_{j,\beta}(r,\theta)
(r\partial_r)^j \partial_\theta^\beta
\end{equation}
where $j$ is an integer, $\beta=(\beta_1,\dots,\beta_{N-1})\in
\N^{N-1}$ is a multi-index and $|\beta|=\sum_{j=1}^{N-1} \beta_j$.

According to \eqref{elliptic-Mazzeo}, the elliptic operator on $X$
corresponding to our operator $\mathcal L_{\A,a}-h$ takes the form
\begin{equation*}
\mathcal L_{\A,a,h}^X:=-r^2
\partial_r^2-(N-1)r\partial_r+L_{\A,a}-r^2 h(r,\theta) .
\end{equation*}
Here by $L_{\A,a}$, we denote the operator on the sphere
$(-i\nabla_{\SN}+\A)^2-a$.

By \cite[Theorem (7.3)]{Mazzeo-1},  if $u$ is a
distributional solution of the equation $\mathcal L_{\A,a,h}^X u=0$
and $r^{-\delta}u(r,\theta)\in L^2(dr d\theta)$, then $u$ admits the
following distributional asymptotic expansion
\begin{equation} \label{M-expansion}
u\sim \sum_{\Re s_j>\delta-\frac12}\sum_{\ell=0}^\infty
\sum_{p=0}^{p_j} r^{s_j+\ell} (\log r)^p u_{j,\ell,p}(\theta)
\end{equation}
where $\{s_j:j\in\Z\setminus\{0\}\}$ coincides with the boundary
spectrum defined in \cite[Definition (2.21)]{Mazzeo-1}. The
numbers $s_j$ are usually called \emph{indicial roots} and in our case
they can be written explicitly in terms of the eigenvalues of the
operator $L_{\A,a}$, i.e.
\begin{equation*}
s_j=-\frac{N-2}2+\text{sign}(j)\sqrt{\left(\frac{N-2}2\right)^2+\mu_{|j|}(\A,a)}
\qquad \text{for all } j\in\Z\setminus\{0\}
\end{equation*}
where $\mu_1(\A,a)\leq\mu_2(\A,a)\leq \mu_3(\A,a)\leq \dots \leq
\mu_k(\A,a)\leq \dots$ denote the eigenvalues of $L_{\A,a}$. For more details
on the meaning of the asymptotic expansion \eqref{M-expansion} see
\cite[Section 7]{Mazzeo-1}.

Let us concentrate our attention on the first term of the
expansion \eqref{M-expansion}, i.e.
\begin{equation} \label{LOG}
r^{s_{j_\delta}} \sum_{p=0}^{p_{j_\delta}} (\log r)^p u_{j_\delta,0,p}(\theta)
\end{equation}
where $j_\delta$ is the smallest value of $j\in \Z$ for which
$s_{j_\delta}>\delta-\frac 12$, see \cite[Theorem (7.3)]{Mazzeo-1}.
This term could be identically zero if $\delta$ is not optimal,
whereas a finer choice of $\delta$ allows selecting the first
nontrivial term in \eqref{M-expansion}.

By using a monotonicity formula approach, in \cite{FFT} we were able to
prove that under the following assumption on  $h$,
i.e.
\begin{equation} \label{hph}
h\in L^{\infty}_{\rm loc}(\Omega\setminus\{0\},\C), \qquad
|h(x)|=O(|x|^{-2+\e})  \quad \mbox{as } |x|\rightarrow 0
\quad\text{for some $\e>0$ ,}
\end{equation}
the presence of logarithmic terms (see \eqref{LOG}) in the leading
part of the asymptotic expansion can be excluded. In the present
paper, we show that the same conclusion can be obtained replacing the 
 pointwise assumption  \eqref{hph} with some  integrability conditions on $h$ and 
its gradient, see
(\ref{hph2})-(\ref{hp:eta_1}).

Here and in \cite{FFT}, the indicial root of the leading term in the
asymptotic expansion of finite energy solutions (namely $H^1$-weak
solutions) to \eqref{u}, is determined by introducing the following
Almgren-type monotonicity function
\begin{equation}\label{eq:37}
  {\mathcal N}_{u,h}(r)=\frac{r\int_{B_r} \big[\big|\nabla u(x)+i
    \frac{A(x/|x|)}{|x|}u(x)\big|^2-\frac{a(x/|x|)}{|x|^2}|u(x)|^2 -(\Re
    h(x))|u(x)|^2 \big] \, dx}{\int_{\partial B_r}|u(x)|^2 \, dS},
  \end{equation}
for any $r\in (0,\overline r)$, with $\overline r\in (0,R)$
sufficiently small. By a blow up argument, we are
able to characterize the indicial root $\gamma$ corresponding to
the leading term in the asymptotic expansion as
\begin{equation}
\gamma=\lim_{r\to 0^+} {\mathcal N}_{u,h}(r) .
\end{equation}
We point out that the monotonicity argument does not need
vanishing of solutions of \eqref{u} outside a small neighborhood
of $r=0$ which is instead required in the Mellin transform
approach used in \cite[Section 7]{Mazzeo-1}. Moreover, here and in
\cite{FFT}, a characterization of the coefficient of the leading
power is given by means of a Cauchy's integral type formula for
$u$, see \eqref{eq:38}.

Let us now describe the integrability type assumptions on
the perturbation $h$ which are required by the forthcoming analysis.
Let $\Omega\subset \R^N$, $N\ge 3$, be a domain containing the origin.
Let $\overline R>0$ be
such that $B_{\overline R}\subset \Omega$ and let $h$ satisfy
\begin{equation} \label{hph2}
h\in L^{\infty}_{\rm loc}(\Omega\setminus\{0\},\C), \qquad \nabla
h\in L^1_{{\rm loc}}(\Omega\setminus\{0\},\C^N) .
\end{equation}
Define, for any $r\in(0,\overline R)$, the two functions
\begin{multline}\label{quotient-1}
\!\!\eta_0(r)=\!\!\sup_{\substack{u\in H^1(B_r)\\u\not\equiv0}}
\frac{\int_{B_r} |h(x)||u(x)|^2 dx} {\int_{B_r} \left|\nabla
u(x)+i\frac{\A(x/|x|)}{|x|}\, u(x)\right|^2 dx- \int_{B_r}
\frac{a(x/|x|)}{|x|^2} \, |u(x)|^2 dx+\frac{N-2}{2r}\int_{\partial
B_r} |u|^2 dS}
\end{multline}

\begin{multline}\label{quotient-2}
  \!\!\eta_1(r)=\!\!\sup_{\substack{u\in H^1(B_r)\\u\not\equiv0}} \frac{\int_{B_r}
    |\Re(x\cdot \nabla h(x))||u(x)|^2 dx} {\int_{B_r} \left|\nabla
      u(x)+i\frac{\A(x/|x|)}{|x|}\, u(x)\right|^2 dx- \int_{B_r}
    \frac{a(x/|x|)}{|x|^2} \, |u(x)|^2 dx+\frac{N-2}{2r}\int_{\partial
      B_r} |u|^2 dS}.
\end{multline}
We observe that, under the assumption
$$
\mu_1(\A,a)>-\left(\frac{N-2}2\right)^2,
$$
the quadratic form appearing at the denominators of the two
quotients in \eqref{quotient-1} and \eqref{quotient-2} is positive
for any $u\in H^1(B_r)\setminus \{0\}$ and for any $r>0$, and its
square root is a norm equivalent to the $H^1(B_r)$-norm (see \cite[Lemma 3.1]{FFT}).

Let us assume that
\begin{equation} \label{hp:eta_0} \lim_{r\to 0^+} \eta_0(r)=0, \qquad
  \frac{\eta_0(r)}{r}\in L^1(0,\overline R), \qquad \frac{1}{r}
  \int_0^r \frac{\eta_0(s)}{s} \, ds\in L^1(0,\overline R) .
\end{equation}
and that
\begin{equation}\label{hp:eta_1}
  \frac{\eta_1(r)}{r}\in
  L^1(0,\overline R), \qquad \frac{1}{r} \int_0^r \frac{\eta_1(s)}{s}
  \, ds\in L^1(0,\overline R) .
\end{equation}
Conditions \eqref{hp:eta_0} and \eqref{hp:eta_1} are satisfied for example if
\begin{equation*}
h\in L^s(B_{\overline R},\C),  \quad |x\cdot \nabla h|\in
L^s(B_{\overline R}), \quad \text{for some } s>N/2
\end{equation*}
or
\begin{equation*}
h\in K_{N,\delta}^{{\rm loc}}(B_{\overline R}) \qquad \text{and} \qquad
\Re(x\cdot \nabla h(x))\in K_{N,\delta}^{{\rm loc}}(B_{\overline R})
\end{equation*}
for some $\delta>0$. Here $K_{N,\delta}^{{\rm loc}}(B_{\overline R})$ denotes a
modified version of the usual Kato class $K_{N}^{{\rm loc}}(B_{\overline R})$
(see \cite{Kato} for the definition of $K_N^{{\rm loc}}(B_{\overline R})$ and
\cite{Kurata} for the definition of $K_{N,\delta}^{{\rm
    loc}}(B_{\overline R})$).

A further aim of the present paper is to point out how the combination
of monotonicity and blow-up techniques provides a powerful tool in the
study of nonlinear problems of the type \eqref{nonlin0}, where $f$ is
a nonlinearity with at most critical growth. In \cite{FFT}, the study
of \eqref{nonlin0} was carried out as follows: a-priori upper bounds
of solutions to \eqref{nonlin0} were first deduced by a classical
iteration scheme, allowing treating the nonlinear term as a linear one
of the type $h(x)u$ with a potential $h$ depending nonlinearly on $u$
but satisfying a suitable pointwise estimate.  The linear result
\cite[Theorem 1.3]{FFT} was thus invoked to prove its nonlinear
version \cite[Theorem 1.6]{FFT}.  In particular, in \cite{FFT} a
nonlinear version of the monotonicity formula was not needed being the
asymptotics for the nonlinear problem deducible from the linear
case. On the other hand, the a-priori pointwise estimate on solutions
of \eqref{nonlin0} needed to reduce the nonlinear problem to a linear
one required the further assumption 
\begin{align} \label{pos-mu10a}
\mu_1(0,a)>-\bigg(\frac{N-2}2\bigg)^{\!\!2},
\end{align}
 see the
statement of \cite[Theorem 1.6]{FFT} and \cite[Theorem 9.4]{FFT}.

In the present paper, we remove condition (\ref{pos-mu10a}) and prove
Theorem \ref{t:asym} below under the less
restrictive positive definiteness condition (\ref{positivity}). Such
improved result is obtained through  a unified approach which allows
 treating simultaneously linear and nonlinear equations.
A similar unified approach was previously introduced in the paper
\cite{FFT2} dealing with  elliptic equations with cylindrical and
many-particle potentials, for which a-priori pointwise estimates
seem to be quite more difficult to be proved, thus requiring a
purely nonlinear approach based on a nonlinear monotonicity
formula.

Let us consider a unified version of \eqref{u} and \eqref{nonlin0},
i.e. an equation of the form
\begin{equation} \label{eq:lin-nonlin}
\mathcal L_{\A,a} u=h(x)\, u+f(x,u), \quad \text{in } \Omega,
\end{equation}
where $h$ satisfies \eqref{hph2}, \eqref{hp:eta_0}, \eqref{hp:eta_1},
 $f$ is of the type
 \begin{equation} \label{special-form} f(x,z)=g(x,|z|^2)z, \qquad
   \text{for a.e. } x\in\Omega, \text{ for all } z\in\C,
\end{equation}
 $g:\Omega\times \R\to \R$
satisfies
\begin{equation}\label{F}
\left\{\!\!
\begin{array}{l}
g\in C^0(\Omega\times [0,+\infty)),\quad G\in C^1(\Omega\times [0,+\infty)),
\\[5pt]
|g(x,s)s|+|\nabla_x G(x,s)\cdot x|\leq C_g(|s|+|s|^{\frac{2^*}{2}}),
\quad\text{for a.e. $x\in\Omega$ and all $s\in\R$},
\end{array}
\right.
\end{equation}
 $G(x,s)=\frac12\int_0^s g(x,t)\,dt$,
$2^*=\frac{2N}{N-2}$ is the critical Sobolev exponent, $C_g>0$ is a
constant independent of $x\in\Omega$ and $s\in\R$, and $\nabla_x G$
denotes the gradient of $G$ with respect to the $x$ variable.

The special form \eqref{special-form} chosen for the function $f$ is
invariant by gauge transformations and hence
very natural in the study of nonlinear Schr\"odinger equations with
magnetic fields, see for example~\cite{el}. We stress that our approach works
 for  very general nonlinearities and also for perturbations of the homogeneous
magnetic potential.

Let us recall the assumptions ${\bf (A.1),(A.2),(A.3),(A.4)}$
already introduced in \cite{FFT}:
\begin{align}
  &{\boldsymbol{\mathcal
      A}}(x)=\dfrac{\A\big(\frac{x}{|x|}\big)}
  {|x|}\quad\text{and}\quad V(x)=-
  \dfrac{a\big(\frac{x}{|x|}\big)}{|x|^2}\tag{\bf A.1}\label{eq:36}
  &\text{\bf\hfil(homogeneity)}\\[5pt]
  &\A\in C^1({\mathbb S}^{N-1},\R^N)\quad\text{and}\quad
  a\in L^{\infty}({\mathbb S}^{N-1},\R)\tag{\bf A.2}\label{eq:reg}
  &\text{\bf (regularity of angular coefficients)}\\[5pt]
&\A(\theta)\cdot\theta=0\quad\text{for all }\theta\in{\mathbb S}^{N-1}.
\tag{\bf A.3}\label{eq:transversality}
  &\text{\bf (transversality)}
\end{align}
\begin{align} \label{positivity}
\quad\!\!\mu_1(\A,a)>-\bigg(\frac{N-2}2\bigg)^{\!\!2}\tag{\bf A.4},
\hskip5.5cm\text{\bf (positive definiteness)}.
\end{align}
An equivalent version of ${\bf (A.4)}$ can be given by introducing
the quantity
\begin{equation} \label{eq:bound}
\Lambda(\A,a):=\sup_{\substack{u\in \Di\\u\not\equiv0}}
\dfrac{{\int_{\R^N}{{|x|^{-2}}{a(x/|x|)}\,|u(x)|^2\,dx}}}
{\int_{\R^N}\big|\nabla u(x)+i
    \frac{\A(x/|x|)}{|x|}u(x)\big|^2\,dx}
\end{equation}
and by taking into account that
\begin{equation} \label{equivalence}
\mu_1(\A,a)>-\bigg(\frac{N-2}2\bigg)^{\!\!2} \quad \text{if and
only if } \quad \Lambda(\A,a)<1 \ ,
\end{equation}
see \cite[Lemma 1.1]{FFT} and \cite[Lemma 2.3]{FFT2}. It is easy
to verity that $\Lambda(\A,a)\geq 0$ and it is zero if and only if
$a\leq 0$ a.e. in $\SN$.

The following theorem characterizes the leading term of the asymptotic
expansion of solutions to (\ref{eq:lin-nonlin}) by means of the
limit of the  associated Almgren-type function
\begin{multline}\label{eq:37-bis}
  {\mathcal N}_{u,h,f}(r)=\frac{r\int_{B_r} \big[\big|\nabla u(x)+i
    \frac{\A(x/|x|)}{|x|}u(x)\big|^2-\frac{a(x/|x|)}{|x|^2}|u(x)|^2\big]
    \, dx}
  {\int_{\partial B_r}|u(x)|^2 \, dS} \\
  -\frac{r\int_{B_r}\big[(\Re h(x))|u(x)|^2+g(x,|u(x)|^2)|u(x)|^2\big]
    \, dx}{\int_{\partial B_r}|u(x)|^2 \, dS} .
 \end{multline}

 \begin{Theorem} \label{t:asym} Let $\Omega\subset\R^N$, $N\geq 3$, be
   a bounded open set containing $0$, ${\bf (A.1),(A.2),(A.3),(A.4)}$
   hold, and $u$ be a weak $H^1(\Omega,\C)$-solution to
   (\ref{eq:lin-nonlin}), $u\not\equiv 0$, with $h$ satisfying
   (\ref{hph2}), (\ref{hp:eta_0}), (\ref{hp:eta_1}) and $f$ satisfying
   (\ref{special-form}) and (\ref{F}) .  Then, letting ${\mathcal
     N}_{u,h,f}(r)$ as in (\ref{eq:37-bis}), there exists $k_0\in \N$,
   $k_0\geq 1$, such that
  \begin{align}\label{eq:20}
    \gamma:=\lim_{r\to 0^+}{\mathcal
      N}_{u,h,f}(r)=-\frac{N-2}2+\sqrt{\bigg(\frac{N-2}{2}
      \bigg)^{\!\!2}+\mu_{k_0}(\A,a)}.
  \end{align}
  Furthermore, if $m\geq
  1$ is the multiplicity of the eigenvalue $\mu_{k_0}(\A,a)$, and
  $\{\psi_i:\,j_0\leq i\leq j_0+m-1\}$ ($j_0\leq k_0\leq j_0+m-1$) is
  an $L^2({\mathbb S}^{N-1},\C)$-orthonormal basis for the eigenspace
  of the operator $L_{\A,a}$ associated to $\mu_{k_0}(\A,a)$, then
\begin{equation} \label{estu}
\lambda^{-\gamma} u(\lambda\theta)\longrightarrow\sum_{i=j_0}^{j_0+m-1}
  \beta_i\psi_{i}(\theta)\quad \text{in }C^{1,\tau}({\mathbb S}^{N-1},\C)
\quad\text{as }\lambda\to 0^+,
\end{equation}
 and
\begin{equation} \label{estgradu}
\lambda^{1-\gamma}\nabla u(\lambda\theta)\longrightarrow
\sum_{i=j_0}^{j_0+m-1} \beta_i
\big(\gamma \psi_{i}(\theta)\theta+\nabla_{{\mathbb S}^{N-1}}
\psi_{i}(\theta)\big)\quad \text{in }C^{0,\tau}({\mathbb S}^{N-1},\C^N)
\quad\text{as }\lambda\to 0^+,
\end{equation}
for any $\tau\in(0,1)$, where
\begin{align}\label{eq:38}
  \beta_i
=\! \int_{{\mathbb S}^{N-1}}\!\bigg[ \frac{u(R\theta)}{R^{\gamma}}+\!
  \int_{0}^R\frac{
    \big(h(s\,\theta)+g(s\theta,|u(s\theta)|^2)\big)u(s\,\theta)}{2\gamma+N-2}
  \bigg(s^{1-\gamma}-\frac{s^{\gamma+N-1}}{R^{2\gamma+N-2}}\bigg)ds
  \bigg]\overline{\psi_{i}(\theta)}\,dS(\theta),
 \end{align}
 for all $R>0$ such that $\overline B_{R}=
\{x\in\R^N:|x|\leq R\}\subset\Omega$
and $(\beta_{j_0},\beta_{j_0+1},\dots,\beta_{j_0+m-1})\neq(0,0,\dots,0)$.
\end{Theorem}

It is worth pointing out how convergence (\ref{estu})
excludes the presence of logarithmic factors in the leading term of
the expansion \eqref{M-expansion}.

Although the proof of Theorem \ref{t:asym} follows essentially the
scheme of Theorem 1.3 in \cite{FFT}, the addition of the nonlinear
term in the Almgren-type function~(\ref{eq:37-bis}) and the
replacement of pointwise assumptions on $h$ with the integral type
ones (\ref{hp:eta_0}--\ref{hp:eta_1}), require some significant
adaptations which are emphasized in Section \ref{Proof-t:asym}. As a
relevant byproduct of Theorem \ref{t:asym} we also obtain the
following pointwise estimate on solutions to (\ref{eq:lin-nonlin}):

\begin{Corollary}\label{cor:stima}
Let  $u$ be a weak $H^1(\Omega,\C)$-solution to
(\ref{eq:lin-nonlin}) and  all the assumptions of Theorem
  \ref{t:asym} hold. Then for any $\Omega'\Subset\Omega$, there exists a
  constant $C=C(\Omega',u)$ such that
\begin{equation} \label{|u|<gamma}
|u(x)|\leq C|x|^{\gamma}  \quad \text{for a.e. every } x\in \Omega'
\end{equation}
 where $\gamma$ is the number defined (\ref{eq:20}).
\end{Corollary}

We point out that Corollary \ref{cor:stima} is a direct consequence of
Theorem \ref{t:asym} which is proved by monotonicity and blow-up
methods, and hence does not require any iterative Brezis-Kato scheme;
in particular, here we can drop the strongest positivity condition
(\ref{pos-mu10a}), which was instead needed in \cite{FFT} to start the
iteration procedure.

\medskip\noindent
{\bf Acknowledgments.}  The authors are indebted to  Frank  Pacard for 
fruitful discussions and for drawing 
to their attention references \cite{LoMc}, \cite{Mazzeo-1},
\cite{Mazzeo-2},  and \cite{pacard}.

\section{Proof of Theorem \ref{t:asym}} \label{Proof-t:asym}

\noindent Solutions to (\ref{eq:lin-nonlin}) satisfy the following Pohozaev-type
identity.

\begin{Proposition} \label{t:pohozaev} Let $\Omega\subset\R^N$, $N\geq 3$,
  be a bounded open set such that  $0\in\Omega$.
Let  $a,\A$ satisfy ${\bf (A.2)}$, and $u$ be a
  weak $H^1(\Omega,\C)$-solution to (\ref{eq:lin-nonlin}) in $\Omega$, with $h$
  satisfying (\ref{hph2}), (\ref{hp:eta_0}--\ref{hp:eta_1}), and $f$
as in (\ref{special-form}--\ref{F}).  Then
\begin{align}\label{eq:26}
-\frac{N-2}2 & \int_{B_r}\bigg[ \bigg|\bigg(\nabla+i\,\frac
  {\A\big({x}/{|x|}\big)}{|x|}\bigg)u\bigg|^2
  -\frac{a\big({x}/{|x|}\big)}{|x|^2}|u|^2\bigg]\,dx\\
  &
  \notag
  \qquad \qquad
+\frac{r}{2}\int_{\partial B_r}\bigg[ \bigg|\bigg(\nabla+i\,\frac
  {\A\big({x}/{|x|}\big)}{|x|}\bigg)u\bigg|^2
  -\frac{a\big({x}/{|x|}\big)}{|x|^2}|u|^2\bigg]\,dS\\
  &
  \notag
=r\int_{\partial B_r}\bigg|\frac{\partial u}{\partial
\nu}\bigg|^2\,dS-\frac 12 \int_{B_r}\Re\big(\nabla h(x) \cdot x
\big) |u(x)|^2\,dx \\
&
  \notag
  \qquad \qquad
-\frac N2\int_{B_r}\Re\big(h(x)\big) |u(x)|^2\,dx+\frac
r2\int_{\partial B_r}\Re\big(h(x)\big) |u(x)|^2\,dS \\
&
\notag
+r\int_{\partial B_r} G(x,|u(x)|^2)\, dS-\int_{B_r} \big(
\nabla_x G(x,|u(x)|^2)\cdot x+NG(x,|u(x)|^2)\big)\, dx
\end{align}
for all $r>0$ such that $\overline{B_{r}}= \{x\in\R^N:|x|\leq
r\}\subset\Omega$, where $\nu=\nu(x)$ is the unit outer normal
vector $\nu(x)=\frac{x}{|x|}$.
\end{Proposition}

\begin{pf} One can proceed  similarly to the proof Theorem 4.1 in \cite{FFT}
by fixing $r\in (0,\overline R)$ and finding a sequence
$\{\delta_n\}\subset (0,r)$ such that $\lim_{n\to +\infty}
\delta_n=0$ and
\begin{equation}\label{eq:33}
\delta_n \int_{\partial B_{\delta_n}}
\bigg[\bigg|\bigg(\nabla+i\,\frac
  {\A\big({x}/{|x|}\big)}{|x|}\bigg)u\bigg|^2 +\frac{|u|^2}{|x|^2}
  +\bigg|\frac{\partial u}{\partial
    \nu}\bigg|^2+\Re\big(h(x)\big) |u(x)|^2
+|G(x,|u(x)|^2)|\bigg]\,dS
    \to 0
\end{equation}
as $n\to+\infty$.
This is possible by the fact that
$\Re\big(h(x)\big) |u(x)|^2,G(x,|u(x)|^2)\in L^1(B_r)$ in view of
(\ref{quotient-1}), (\ref{hp:eta_0}) and (\ref{F}).

 By ${\bf (A.2)}$ and \eqref{hph2} we deduce
that $u\in C^{1,\tau}_{\rm loc}(\Omega\setminus\{0\},\C)$ for any
$\tau\in(0,1)$ and $h\in W^{1,1}_{\rm
loc}(\Omega\setminus\{0\},\C)$ and hence, integrating by parts, we
obtain
\begin{align*}
\int_{B_r\setminus B_{\delta_n}}& \Re\big(h(x) u(x) (x\cdot
\overline{\nabla u(x)}) \big) \, dx\\
&= -\frac 12 \int_{B_r\setminus B_{\delta_n}}\Re\big(\nabla h(x)
\cdot x \big) |u(x)|^2\,dx -\frac N2\int_{B_r\setminus
B_{\delta_n}}\Re\big(h(x)\big) |u(x)|^2\,dx
\\
& \quad +\frac r2\int_{\partial B_r}\Re\big(h(x)\big) |u(x)|^2\,dS
-\frac{\delta_n}2\int_{\partial B_{\delta_n}}\Re\big(h(x)\big)
|u(x)|^2\,dS .
\end{align*}
Passing to the limit as $n\to +\infty$, by \eqref{hp:eta_0},
\eqref{hp:eta_1}, and \eqref{eq:33} we obtain
\begin{align*}
\lim_{n\to +\infty}&\int_{B_r\setminus B_{\delta_n}}
\Re\big(h(x) u (x\cdot \overline{\nabla u(x)}) \big) \, dx\\
&=-\frac 12 \int_{B_r}\Re\big(\nabla h(x) \cdot x \big)
|u(x)|^2\,dx -\frac N2\int_{B_r}\Re\big(h(x)\big) |u(x)|^2\,dx
+\frac r2\int_{\partial B_r}\Re\big(h(x)\big) |u(x)|^2\,dS
 .
\end{align*}
The proof of the proposition then follows proceeding as in the
proof of Theorem 4.1 in \cite{FFT} and Proposition A.1 in \cite{FFT2}.
\end{pf}

Proceeding as in \cite{FFT}, one can show that, under the assumptions
${\bf (A.2)}, {\bf (A.3)}, {\bf (A.4)}$, and
 \eqref{hp:eta_0}, there exists $\overline
r\in (0,\overline R)$ such that the function $H(r)=r^{1-N}
\int_{\partial B_r} |u|^2 dS$ is strictly positive for any $r\in
(0,\overline r)$ and $\sup_{r\in(0,\overline r)}\eta_0(r)<+\infty$. In
this way, if $D$ is the function defined by
\begin{multline*}
  D(r)=\frac{1}{r^{N-2}}\int_{B_r} \left[\left|\nabla u(x)+i
      \frac{\A(x/|x|)}{|x|}u(x)\right|^2-
    \frac{a(x/|x|)}{|x|^2}|u(x)|^2\right] \, dx \\
  -\frac{1}{r^{N-2}} \int_{B_r}\big[(\Re
  h(x))|u(x)|^2+g(x,|u(x)|^2)|u(x)|^2\big] \, dx,
 \end{multline*}
then the quotient
\begin{equation} \label{well-def}
\mathcal N(r):=\mathcal N_{u,h,f}(r)=\frac{D(r)}{H(r)}, \quad
\text{for a.e. } r\in(0,\overline r),
\end{equation}
 is well defined.
Arguing as in \cite[(52)]{FFT}, it is easy to verify that
\begin{equation}\label{eq:DprimoH}
D(r)=\frac{r}2H'(r)\quad \text{for a.e. } r\in(0,\overline r).
\end{equation}

\begin{Lemma} \label{mono} Let $\Omega\subset\R^N$, $N\geq 3$, be a
  bounded open set such that $0\in\Omega$, $a,\A$ satisfy ${\bf (A.2),
    (A.3)}$, ${\bf (A.4)}$, and $u\not\equiv 0$ be a weak
  $H^1(\Omega,\C)$-solution to (\ref{eq:lin-nonlin}) in $\Omega$, with
  $h$ satisfying (\ref{hph2}), (\ref{hp:eta_0}--\ref{hp:eta_1}), and
  $f$ satisfying (\ref{special-form}--\ref{F}).  Then, letting
  ${\mathcal N}$ as in (\ref{well-def}), there holds ${\mathcal N}\in
  W^{1,1}_{{\rm loc}}(0,\overline r)$ and
\begin{gather} \label{formulona}
{\mathcal N'}(r)=\nu_1(r)+\nu_2(r)
\end{gather}
in a distributional sense and for a.e. $r\in (0,\overline r)$,
where
\begin{equation} \label{eq:nu1}
\nu_1(r)= \frac{2r\left[
      \left(\int_{\partial B_r} \left|\frac{\partial
            u}{\partial\nu}\right|^2 dS\right) \cdot
      \left(\int_{\partial B_r} |u|^2 dS\right)-\left(\int_{\partial
          B_r} \Re\left(u\frac{\partial \overline u}{\partial
            \nu}\right) dS\right)^{\!2} \right]}
  {\left(\int_{\partial B_r} |u|^2 dS\right)^2}
\end{equation}
and
\begin{align} \label{eq:nu2}
\nu_2(r)=&-\frac{ \int_{B_r}\Re\big(2h(x)+\nabla h(x) \cdot x
\big) |u(x)|^2\,dx } {\int_{\partial B_r} |u|^2 dS}\\
\notag &+\frac{r\int_{\partial
B_r}\big(2G(x,|u(x)|^2)-g(x,|u(x)|^2)|u(x)|^2\big)\, dS}
{\int_{\partial B_r}|u|^2 \, dS} \\[10pt]
\notag &
+\frac{\int_{B_r}
\big((N-2)g(x,|u(x)|^2)|u(x)|^2-2NG(x,|u(x)|^2) -2\nabla_x
G(x,|u(x)|^2)\cdot x\big)\, dx} {\int_{\partial B_r}|u|^2 \, dS}.
\end{align}
\end{Lemma}

\begin{pf}
One can proceed exactly as in the proof of Lemma 5.4 in \cite{FFT}
by using the Pohozaev-type identity \eqref{eq:26} in place of (32)
in \cite{FFT}.
\end{pf}

The following proposition provides an a-priori super-critical
summability of solutions to (\ref{eq:lin-nonlin}) which will allow
including the critical growth case in the Almgren type monotonicity
formula.

\begin{Proposition} \label{SMETS}
Let $\Omega\subset \R^N$, $N\geq 3$ be a bounded open set such
that $0\in \Omega$, $a$, $\A$ satisfy ${\bf (A.2)}, {\bf (A.3)},
{\bf  (A.4)}$,
and $u$ be a $H^1(\Omega,\C)$-weak solution to
\begin{equation} \label{eq:ahV}
\mathcal L_{\A,a} u(x)=h(x)u(x)+V(x)u(x), \quad \text{in } \Omega,
\end{equation}
with $h$ satisfying (\ref{hph2}), (\ref{hp:eta_0}--\ref{hp:eta_1})
 and $V\in L^{N/2}(\Omega,\C)$. Letting
$$
q_{\rm lim}:=
\begin{cases}
  \frac{2^*}{2} \min\left\{\frac{4}{\Lambda(\A,a)}-2,2^*\right\}, &
  \text{if }\Lambda(\A,a)>0, \\[5pt]
  \frac{(2^*)^2}2, & \text{if } \Lambda(\A,a)=0 ,
\end{cases}
$$
then for any $1\leq q<q_{\rm lim}$ there exists $r_q>0$, depending
only on $N,\A,a,q,h$ such that $B_{r_q}\subset \Omega$ and $u\in
L^q(B_{r_q},\C)$.
\end{Proposition}

\begin{pf}
By ${\bf (A.4)}$ and \eqref{equivalence} we have that $\frac{2}{2^*}
q_{\rm lim}>2$. For any $2<\tau<\frac{2}{2^*} q_{\rm lim}$, define
$C(\tau):=\frac{4}{\tau+2}$ and let $\ell_\tau>0$ be so large
 that
\begin{equation} \label{elleq}
\bigg(\int\limits_{|V(x)| \geq\ell_\tau}\!\!\!\!\!\!
|V(x)|^{\frac{N}{2}} \, dx
\bigg)^{\!\!\frac{2}{N}}<\frac{S(\A)(C(\tau)-\Lambda(\A,a))}2
\end{equation}
where
$$
S(\A):=\inf_{\substack{v\in \Di \\ v\not\equiv 0 }}
\frac{\int_{\R^N} \left|\nabla
v(x)+i\frac{\A(x/|x|)}{|x|}v(x)\right|^2 \, dx}{\left(\int_{\R^N}
|v(x)|^{2^*} \, dx\right)^{\frac{2}{2^*}}}>0 .
$$
Let $r>0$ be such that $B_r\subset \Omega$.  For any $w\in
H^1_0(B_{r},\C)$, by H\"older and Sobolev inequalities and
\eqref{elleq}, we have
\begin{align} \label{estV}
\int_{B_{r}} & |V(x)||w(x)|^2 \, dx= \int \limits_{B_r\cap
\{|V(x)| \leq\ell_\tau\}}\!\!\!\!\!\!\ |V(x)||w(x)|^2 \, dx + \int
\limits_{B_r\cap \{|V(x)| \geq\ell_\tau\}}\!\!\!\!\!\!\
|V(x)||w(x)|^2 \, dx
\\
\notag & \leq \ell_\tau\int_{B_r} |w(x)|^2 \, dx +
\bigg(\int\limits_{|V(x)| \geq\ell_\tau}\!\!\!\!\!\!
|V(x)|^{\frac{N}{2}} \, dx \bigg)^{\!\!\frac{2}{N}} \left(\int_{B_r}
|w(x)|^{2^*}\, dx\right)^{\!\!\frac{2}{2^*}}
\\
\notag & \leq \ell_\tau \int_{B_r} |w(x)|^2 \,
dx+\frac{C(\tau)-\Lambda(\A,a)}2 \int_{B_r} \left|\nabla
w(x)+i\frac{\A(x/|x|)}{|x|}w(x)\right|^2 \, dx .
\end{align}
Let $\rho\in C^\infty_c (B_r,\R)$ be such that $\rho\equiv 1$ in
$B_{r/2}$ and define $v(x):=\rho(x)u(x)\in H^1_0(B_r,\C)$. Then
$v$ is a $H^1(\Omega,\C)$-weak solution of the equation
\begin{equation} \label{eq:g}
\mathcal L_{\A,a} v(x)=h(x)v(x)+V(x)v(x)+g(x) \qquad \text{in }
\Omega
\end{equation}
where $g(x)=-u(x)\Delta \rho(x)-2\nabla u(x)\cdot
\nabla\rho(x)-2iu(x)\frac{\A(x/|x|)}{|x|}\cdot \nabla \rho(x)\in
L^2(B_r,\C)$.
For any $n\in\N$, $n\geq 1$, let us define the function
$v^n:=\min\{|v|,n\}$. Testing \eqref{eq:g} with $(v^n)^{\tau-2}
\overline v\in H^1_0(B_r,\C)$ we obtain
\begin{align} \label{vn}
\int_{B_r} & \!\!\!(v^n(x))^{\tau-2} \left|\nabla
v(x)+i{\textstyle{\frac{\A(x/|x|)}{|x|}}}v(x) \right|^2 \, dx +(\tau-2)
\int_{B_r}\!\!\! (v^n(x))^{\tau-2} |\nabla |v(x)||^2
\chi_{\{|v(x)|<n\}}(x) \, dx \\
\notag & \qquad -\int_{B_r} \frac{a(\frac{x}{|x|})}{|x|^2}
(v^n(x))^{\tau-2} |v(x)|^2\,
dx \\
\notag & =\int_{B_r} \Re(h(x))(v^n(x))^{\tau-2} |v(x)|^2\, dx
+\int_{B_r}\Re( V(x))(v^n(x))^{\tau-2} |v(x)|^2\, dx\\
\notag & \qquad + \int_{B_r} \Re(g(x)(v^n(x))^{\tau-2} \overline
v(x))\, dx .
\end{align}
Since
\begin{multline*}
\left|\nabla((v^n)^{\frac{\tau}{2}-1}v)+i\frac{\A(x/|x|)}{|x|}(v^n)^{\frac{\tau}{2}-1}v
\right|^2 \\
= (v^n)^{\tau-2} \left|\nabla
v+i\frac{\A(x/|x|)}{|x|}v\right|^2+\frac{(\tau-2)(\tau+2)}4
(v^n)^{\tau-2}  |\nabla |v||^2 \chi_{\{|v(x)|<n\}},
\end{multline*}
then by \eqref{vn}, \eqref{eq:bound}, (\ref{quotient-1}), and
\eqref{estV} with $w=(v^n)^{\frac{\tau}{2}-1}v$, we obtain for any
$r>0$ small enough such that $\eta_0(r)<1$,
\begin{gather} \label{vn2}
C(\tau) \int_{B_r}
  \left|\nabla((v^n)^{\frac{\tau}{2}-1}v)+
    i\frac{\A(x/|x|)}{|x|}(v^n)^{\frac{\tau}{2}-1}v
  \right|^2\, dx \\
  \notag \leq \int_{B_r} \frac{a(\frac{x}{|x|})}{|x|^2}
  |(v^n(x))^{\frac{\tau}{2}-1}v(x)|^2 \, dx+ \int_{B_r}
  \Re(h(x))|(v^n(x))^{\frac{\tau}{2}-1}v(x)|^2\,
  dx  \\
  \notag \quad +\int_{B_r}
  \Re(V(x))|(v^n(x))^{\frac{\tau}{2}-1}v(x)|^2\,
  dx + \int_{B_r} \Re\big(g(x)(v^n(x))^{\tau-2} \overline{v(x)}\big)\, dx \\
  \notag \leq \bigg[\Lambda(\A,a)(1-\eta_0(r))+\eta_0(r)
  +\frac{C(\tau)-\Lambda(\A,a)}2\bigg] \int_{B_r}
  \left|\nabla((v^n)^{\frac{\tau}{2}-1}v)+
    i\frac{\A(x/|x|)}{|x|}(v^n)^{\frac{\tau}{2}-1}v \right|^2 dx
  \\
  \notag \quad +\ell_\tau \int_{B_r} (v^n(x))^{\tau-2} |v(x)|^2 \,dx
  +\int_{B_r} |g(x)|(v^n(x))^{\tau-2} |v(x)|\, dx .
\end{gather}
Let us consider the last term in the right hand side of
\eqref{vn2}. Since $g\in L^2(B_r,\C)$, then by H\"older inequality
\begin{align*}
\int_{B_r} & |g(x)|(v^n(x))^{\tau-2} |v(x)|\, dx \leq
\|g\|_{L^2(\Omega,\C)} \left(\int_{B_r} (v^n(x))^{2\tau-4}|v(x)|^2
\,
dx\right)^{\!\!\frac{1}2}\\
\notag & =\|g\|_{L^2(\Omega,\C)} \left(\int_{B_r}
(v^n(x))^{\frac{2(\tau-1)(\tau-2)}\tau}
(v^n(x))^{\frac{2(\tau-2)}\tau}|v(x)|^2
\, dx\right)^{\!\!\frac{1}2} \\
\notag &\leq \|g\|_{L^2(\Omega,\C)} \left(\int_{B_r}
|(v^n(x))^{\frac{\tau}2-1} v(x)|^{\frac{4(\tau-1)}\tau} \,
dx\right)^{\!\!\frac{1}2}
\end{align*}
and, since $\frac{4(\tau-1)}\tau<2^*$ for any $\tau<\frac{2}{2^*}
q_{\rm lim}$, by H\"older inequality, Sobolev embedding, and Young
inequality, we obtain
\begin{align} \label{est-g}
  &\int_{B_r}  |g(x)|(v^n(x))^{\tau-2} |v(x)|\, dx\\
  \notag & \leq \|g\|_{L^2(\Omega,\C)}
  \bigg(\frac{\omega_{N-1}}N\bigg)^{\!\!\frac{1}2-\frac{2(\tau-1)}{2^*\tau}}
  r^{\frac{N}2-\frac{2N(\tau-1)}{2^*\tau}} \left(\int_{B_r}
    |(v^n(x))^{\frac{\tau}2-1} v(x)|^{2^*} \,
    dx\right)^{\!\!\frac{2(\tau-1)}{2^*\tau}} \\
  \notag & \leq \|g\|_{L^2(\Omega,\C)}
  \bigg(\frac{\omega_{N-1}}N\bigg)^{\!\!\frac{1}2-\frac{2(\tau-1)}{2^*\tau}}
  r^{\frac{N}2-\frac{(N-2)(\tau-1)}{\tau}}
  S(\A)^{-\frac{\tau-1}\tau} \times \\
  \notag & \qquad \times \left(\int_{B_r}
    \left|\nabla((v^n)^{\frac{\tau}{2}-1}v)+
      i\frac{\A(x/|x|)}{|x|}(v^n)^{\frac{\tau}{2}-1}v
    \right|^2 dx\right)^{\!\!\frac{\tau-1}{\tau}}  \\
  \notag & \leq \frac{\tau-1}\tau
  \bigg(\frac{\omega_{N-1}}N\bigg)^{\frac{\tau}{2(\tau-1)}-\frac{2}{2^*}}
  \frac{r^{\frac{N\tau}{2(\tau-1)}-N+2}}{S(\A)} \int_{B_r}
  \left|\nabla((v^n)^{\frac{\tau}{2}-1}v)+
    i\frac{\A(x/|x|)}{|x|}(v^n)^{\frac{\tau}{2}-1}v \right|^2
  dx  \\
  \notag & \qquad +\frac{1}{\tau} \|g\|^\tau_{L^2(\Omega,\C)},
\end{align}
where $\omega_{N-1}$ denotes the volume of the unit sphere ${\mathbb
  S}^{N-1}$, i.e. $\omega_{N-1}=\int_{{\mathbb S}^{N-1}}dS(\theta)$.
Inserting \eqref{est-g} into \eqref{vn2} we obtain
\begin{align*}
\bigg[& \frac{C(\tau)-\Lambda(\A,a)}2-\eta_0(r)-\frac{\tau-1}\tau
\bigg(\frac{\omega_{N-1}}N\bigg)^{\!\!\frac{\tau}{2(\tau-1)}-\frac{2}{2^*}}
r^{\frac{N\tau}{2(\tau-1)}-N+2}  S(\A)^{-1} \bigg]\times
\\ \notag &
\qquad\qquad\qquad \times \int_{B_r}
\left|\nabla((v^n)^{\frac{\tau}{2}-1}v)+i\frac{\A(x/|x|)}{|x|}(v^n)^{\frac{\tau}{2}-1}v
\right|^2
 dx \\
\notag & \leq \frac{1}{\tau} \|g\|^\tau_{L^2(\Omega,\C)}
+\ell_\tau \int_{B_r} (v^n(x))^{\tau-2} |v(x)|^2 \,dx
\end{align*}
and, by Sobolev embedding,
\begin{align} \label{lim-q}
S(\A) \bigg[&
\frac{C(\tau)-\Lambda(\A,a)}2-\eta_0(r)-\frac{\tau-1}\tau
\bigg(\frac{\omega_{N-1}}N\bigg)^{\!\!\frac{\tau}{2(\tau-1)}-\frac{2}{2^*}}
r^{\frac{N\tau}{2(\tau-1)}-N+2}  S(\A)^{-1} \bigg]\times \\
\notag & \qquad\times
 \left(
\int_{B_r}
(v^n(x))^{\frac{2^*}2\tau-2^*} |v(x)|^{2^*} \, dx\right)^{\!\!2/2^*} \\
\notag & \leq \frac{1}{\tau} \|g\|^\tau_{L^2(\Omega,\C)}
+\ell_\tau \int_{B_r} (v^n(x))^{\tau-2} |v(x)|^2 \,dx.
\end{align}
 Since $\tau<\frac{2}{2^*}
q_{\rm lim}$ then $C(\tau)-\Lambda(\A,a)$ is positive and
$\frac{N\tau}{2(\tau-1)}-N+2$ is also positive. Moreover by
\eqref{hp:eta_0}, $\lim_{r\to 0^+} \eta_0(r)=0$.
 Hence we
may fix $r$ small enough in such a way that the left hand side of
\eqref{lim-q} becomes positive. Since $v\in L^\tau(B_{r},\C)$,
letting $n\to +\infty$, the right hand side of \eqref{lim-q}
remains bounded and hence, by Fatou Lemma, we infer that $v\in
L^{\frac{2^*}2 \tau}(B_{r},\C)$. Since $\rho\equiv 1$ in $B_{r/2}$,
we may conclude that $u\in L^{\frac{2^*}2 \tau}(B_{r/2},\C)$. This
completes the proof of the lemma.
\end{pf}

\noindent
According to the previous proposition, we may fix from now on
a weak $H^1$-solution $u$  to (\ref{eq:lin-nonlin}),
$$
2^*<q<q_{\rm lim},
$$
and $r_q$ in such a way that $u\in L^q(B_{r_q})$.  We omit the proof
of the following lemma which can be deduced in a quite standard way by
combining Hardy-Sobolev inequalities with boundary terms (see \cite[\S
3]{FFT}) with assumptions (\ref{hp:eta_0}) and (\ref{F}).

\begin{Lemma} \label{l:stimasotto}
Under the same assumptions as in Lemma \ref{mono}, there exist
$\tilde r\in (0,\min\{\overline r,r_q\})$ and a positive constant
$\overline{C}=\overline{C}(N,\A,a,h,f,u)>0$ depending on $N$,
$\A$, $a$, $h$, $f$, $u$ but independent of $r$ such that
\begin{multline}\label{eq:47}
  \int_{B_r}  \left[\left|\nabla u(x)+i
      \frac{\A(x/|x|)}{|x|}u(x)\right|^2-
    \frac{a(x/|x|)}{|x|^2}|u(x)|^2\right] \, dx \\
   -\int_{B_r}\big[(\Re
  h(x))|u(x)|^2+g(x,|u(x)|^2)|u(x)|^2\big] \, dx
  \\[10pt]
  \geq  -\frac{N-2}{2r} \int_{\partial B_r} |u(x)|^2
  dS+\overline{C}
  \bigg(\int_{B_r}|u(x)|^{2^*}dx\bigg)^{\!\!\frac2{2^*}}
  \\[10pt]
   +\overline{C} \bigg( \int_{B_r} \left[\left|\nabla u(x)+i
      \frac{\A(x/|x|)}{|x|}u(x)\right|^2-\frac{a(x/|x|)}{|x|^2}|u(x)|^2\right]
  \, dx +\frac{N-2}{2r} \int_{\partial B_r} |u(x)|^2 dS\bigg)
\end{multline}
and
\begin{equation}\label{Nbelow}
   {\mathcal N}(r)>-\frac{N-2}{2}
 \end{equation}
for every $r\in(0,\tilde r)$.
\end{Lemma}

\noindent The term $\nu_2$ introduced in Lemma \ref{mono} can be estimated as follows.
\begin{Lemma} \label{l:stima_nu2} Under the same assumptions as in
  Lemma \ref{mono}, let $\tilde r $ be as in Lemma \ref{l:stimasotto} and
  $\nu_2$ as in (\ref{eq:nu2}). Then there exist a positive constant
  $C_1>0$ depending on $N,q,C_g,\overline{C},\tilde
  r,\|u\|_{L^q(B_{\tilde r},\C)}$ and a function $\omega\in L^1(0,\tilde r)$,
  $\omega\geq 0$ a.e. in $(0,\tilde r)$, such that
$$
|\nu_2(r)|\leq C_1\left[{\mathcal
      N}(r)+\frac{N}{2}\right]\Big[r^{-1}(\eta_0(r)+\eta_1(r))+r^{-1+
\frac{2(q-2^*)}{q}}+\omega(r)\Big]
$$
for a.e. $r\in (0,\tilde r)$ and
$$
\int_0^r \omega(s)\,ds\leq \frac
{\|u\|_{L^{2^*}(\Omega)}^{2^*(1-\alpha)}}{1-\alpha} \,
r^{\frac{N(q-2^*)}q(\alpha-\frac{2}{2^*})}
$$
for all $r\in (0,\tilde r)$ and for some $\alpha$ satisfying
$\frac{2}{2^*}<\alpha<1$.
\end{Lemma}

\begin{pf} The estimates on the terms in \eqref{eq:nu2} involving $g$
  and $G$ can be obtained by using Proposition \ref{SMETS}, Lemma
  \ref{l:stimasotto} and proceeding as in the proof of
  \cite[Lemma 5.6]{FFT2}.

Here we only estimate the term in \eqref{eq:nu2}
which involves the function $h$ and its gradient. From
\eqref{quotient-1}, \eqref{quotient-2} and (\ref{eq:47}) we deduce
that
\begin{align*}
\left|\int_{B_r} (2h(x)+\nabla h(x)\cdot
x)|u(x)|^2\,dx\right|&\leq (2\eta_0(r)+\eta_1(r))\overline
C^{-1}r^{N-2}\left[D(r)+\frac{N-2}2 H(r)\right]
\end{align*}
and, therefore,
\begin{align}\label{B00}
  \left|\frac{\int_{B_r} (2h(x)+\nabla h(x)\cdot x)|u(x)|^2\,dx}
    {\int_{\partial B_r}|u|^2 \, dS}\right|&
\leq \overline C^{-1} r^{-1}(2\eta_0(r)+\eta_1(r)) \left[\mathcal
N(r)+\frac{N-2}2 \right]
\end{align}
for all $r\in (0,\tilde r)$.
\end{pf}

\begin{Lemma} \label{gamma}
Under the same assumptions as in Lemma \ref{mono}, the limit
$$
\gamma:=\lim_{r\rightarrow 0^+} {\mathcal N}(r)
$$
exists and is finite.
\end{Lemma}
\begin{pf}
  From Schwarz's inequality, the function $\nu_1$ defined in
  (\ref{eq:nu1}) is nonnegative. Furthermore, by Lemma
  \ref{l:stima_nu2} and assumptions (\ref{hp:eta_0}) and
  (\ref{hp:eta_1}), $\frac{\nu_2}{{\mathcal N}+N/2}\in L^1(0,\tilde
  r)$. Hence, from (\ref{formulona}) and integration we deduce that
${\mathcal N}$ is bounded in $(0,\tilde r)$, thus implying, in view of Lemma
\ref{l:stima_nu2}, that $\nu_2\in   L^1(0,\tilde
  r)$. Therefore ${\mathcal N}'$ turns out to be an integrable  perturbation
of a nonnegative function and hence ${\mathcal N}(r)$ admits a finite limit as
$r\to 0^+$.
For more details, we refer the reader to
Lemmas 5.7 and 5.8 in~\cite{FFT2}.~\end{pf}

\noindent A first consequence of the convergence of  ${\mathcal
  N}$ at $0$ is the following estimate of $H$ from above.

\begin{Lemma} \label{l:uppb}
  Under the same assumptions as in Lemma \ref{mono}, let
  $\gamma:=\lim_{r\rightarrow 0^+} {\mathcal N}(r)$ be as in Lemma \ref{gamma}.
 Then there exists a constant
$K_1>0$ such that
\begin{equation} \label{1stest}
H(r)\leq K_1 r^{2\gamma}  \quad \text{for all } r\in (0,\bar r).
\end{equation}
\end{Lemma}
\begin{pf}
From (\ref{eq:DprimoH}), (\ref{formulona}), and Schwarz's inequality,
it follows that
\begin{equation*}
\frac{H'(r)}{H(r)}=\frac2r {\mathcal N}(r)\geq
\frac{2\gamma}r+\frac2r\int_0^r\nu_2(s)\,ds.
\end{equation*}
By Lemma \ref{l:stima_nu2}, assumptions
(\ref{hp:eta_0}--\ref{hp:eta_1}), and boundedness of ${\mathcal N}$,
we have that $r\mapsto\frac1r\int_0^r\nu_2\in L^1(0,\tilde r)$. Hence
the conclusion follows from integration.
\end{pf}

\noindent We omit the proof of the following lemma which follows closely
the blow up scheme developed in \cite[Lemma 6.1]{FFT}.
\begin{Lemma}\label{l:blowup}
Under the same assumptions as in Lemma \ref{mono}, the following holds true:
\begin{itemize}
\item[\rm (i)] there exists $k_0\in \N$ such that
  $\gamma=-\frac{N-2}2+\sqrt{\left(\frac{N-2}{2}\right)^2+\mu_{k_0}(\A,a)}$;
\item[\rm (ii)] for every sequence $\lambda_n\to0^+$, there exist
a subsequence $\{\lambda_{n_k}\}_{k\in\N}$ and an eigenfunction
$\psi$ of the operator $L_{\A,a}$ associated to the eigenvalue
$\mu_{k_0}(\A,a)$ such that $\|\psi\|_{L^{2}({\mathbb
S}^{N-1},\C)}=1$ and
\[
\frac{u(\lambda_{n_k}x)}{\sqrt{H(\lambda_{n_k})}}\to
|x|^{\gamma}\psi\Big(\frac x{|x|}\Big)
\]
weakly in $H^1(B_1,\C)$, strongly in $H^1(B_r,\C)$ for every
$0<r<1$, and in $C^{1,\tau}_{\rm loc}(B_1\setminus\{0\},\C)$ for
any $\tau\in (0,1)$.
\end{itemize}
\end{Lemma}

\noindent A first step towards the description of the behavior of $H$
as $r\to 0^+$ is the following lemma, whose proof is similar to
\cite[Lemma 6.6]{FFT2}.
\begin{Lemma} \label{l:limite}
Under the same assumptions as in Lemma
  \ref{mono} and letting $\gamma:=\lim_{r\rightarrow 0^+} {\mathcal
    N}(r)\in \R$ as in Lemma \ref{gamma}, the limit
\[
\lim_{r\to0^+}r^{-2\gamma}H(r)
\]
exists and it is finite.
\end{Lemma}

Under the integral type assumptions
(\ref{hp:eta_0}--\ref{hp:eta_1}), the proof that
$\lim_{r\to0^+}r^{-2\gamma}H(r)>0$ is more delicate than it was
under the pointwise conditions required in \cite{FFT} and a new
argument is needed to prove it.

\begin{Lemma}\label{phi>0}
  Suppose that all the assumptions of Lemma \ref{mono} hold true.  Let
  $k_0$ be as in Lemma \ref{l:blowup} and let $j_0,m\in \N$, $j_0,m\ge
  1$ such that $m$ is the multiplicity of $\mu_{k_0}(\A,a)$, $j_0\leq
  k_0\leq j_0+m-1$ and
  $\mu_{j_0}(\A,a)=\mu_{j_0+1}(\A,a)=\cdots=\mu_{j_0+m-1}(\A,a)=\mu_{k_0}(\A,a)$.
  Let $\{\psi_i:j_0\leq i \leq j_0+m-1\}$ be an
  $L^2(\SN,\C)$-orthonormal basis for the eigenspace of the operator
  $L_{\A,a}$ associated to $\mu_{k_0}(\A,a)$.  Then for any sequence
  $\lambda_n\to 0^+$ there exists $i\in \{j_0,\dots,j_0+m-1\}$ such
  that
  \begin{equation*}
    \liminf_{n\to +\infty} \frac{\left|\int_{\SN} u(\lambda_n \theta)
        \overline{\psi_i(\theta)} \,
        dS(\theta)\right|}{\sqrt{H(\lambda_n)}}>0 .
  \end{equation*}
\end{Lemma}

\begin{pf}Suppose by contradiction that there exists a sequence
  $\lambda_n\to 0^+$ such that
\begin{equation*}
  \liminf_{n\to +\infty} \frac{\left|\int_{\SN} u(\lambda_n \theta)
      \overline{\psi_i(\theta)} \,
      dS(\theta)\right|}{\sqrt{H(\lambda_n)}}=0
\end{equation*}
for all $i\in \{j_0,\dots,j_0+m-1\}$.
By Lemma \ref{l:blowup} we deduce that there exist a subsequence
$\{\lambda_{n_k}\}$ and an eigenfunction $\psi$ of the operator
$L_{\A,a}$ corresponding to the eigenvalue $\mu_{k_0}(\A,a)$ with
$\|\psi\|_{L^2(\SN,\C)}=1$, such that
\[
\frac{u(\lambda_{n_k}\theta)}{\sqrt{H(\lambda_{n_k})}}\to
\psi(\theta)
\]
strongly in $L^2(\SN)$ and
\begin{equation*}
\lim_{k\to +\infty} \int_{\SN} \frac{u(\lambda_{n_k}
\theta)}{\sqrt{H(\lambda_{n_k})}}\, \overline{\psi_i(\theta)} \,
dS(\theta)=0 .
\end{equation*}
Therefore
\begin{equation} \label{nullo}
 \int_{\SN} \psi(\theta)\overline{\psi_i(\theta)} \,
dS(\theta)=\lim_{k\to +\infty} \int_{\SN} \frac{u(\lambda_{n_k}
\theta)}{\sqrt{H(\lambda_{n_k})}} \, \overline{\psi_i(\theta)} \,
dS(\theta)=0
\end{equation}
for any $i\in \{j_0,\dots,j_0+m-1\}$. Hence  $\psi\equiv 0$, thus giving
rise to a contradiction.
\end{pf}

\begin{Lemma} \label{l:limitepositivo}
Under the same assumptions as in Lemma
  \ref{mono} and letting $\gamma:=\lim_{r\rightarrow 0^+} {\mathcal
    N}(r)\in \R$ as in Lemma \ref{gamma},  there holds
\[
\lim_{r\to0^+}r^{-2\gamma}H(r)>0.
\]
\end{Lemma}
\begin{pf}
  For the sake of completeness, we report here part of the proof of
  Lemma 6.5 in \cite{FFT}.  Let $0<R<\frac{\tilde r}2$, $\tilde r$ as in
  Lemma \ref{l:stimasotto}, and, for any $k\in\N\setminus\{0\}$, let $\psi_k$
  be a $L^2$-normalized eigenfunction of the operator $L_{\A,a}$ on
  the sphere associated to the $k$-th eigenvalue $\mu_{k}(\A,a)$, i.e.
  satisfying
\begin{equation}\label{eq:2rad}
\begin{cases}
L_{\A,a}\psi_k(\theta)
=\mu_k(\A,a)\,\psi_k(\theta),&\text{in }{\mathbb S}^{N-1},\\[3pt]
\int_{{\mathbb S}^{N-1}}|\psi_k(\theta)|^2\,dS(\theta)=1.
\end{cases}
\end{equation}
We can choose the functions $\psi_k$ in such a way that they form
an orthonormal basis of $L^2({\mathbb S}^{N-1},\C)$, hence $u$ and
$hu+g(x,|u|^2)u$
 can be expanded as
\begin{gather} \label{expansion}
  u(x)=u(\lambda\,\theta)=\sum_{k=1}^\infty\varphi_k(\lambda)\psi_k(\theta),\\
   h(x)u(x)+g(x,|u(x)|^2)u(x)=h(\lambda\,\theta)u(\lambda\,\theta)
  +g(\lambda\,\theta,|u(\lambda\,\theta)|^2)u(\lambda\,\theta)
  =\sum_{k=1}^\infty\zeta_k(\lambda)\psi_k(\theta),
\end{gather}
where $\lambda=|x|\in(0,R]$, $\theta=x/|x|\in{{\mathbb S}^{N-1}}$, and
\begin{equation}\label{eq:22}
  \varphi_k(\lambda)=\!\int_{{\mathbb S}^{N-1}}\!\!u(\lambda\,\theta)
  \overline{\psi_k(\theta)}\,dS(\theta),
  \quad
  \zeta_k(\lambda)=\!\int_{{\mathbb S}^{N-1}}\!\!
  \big(h(\lambda\,\theta) +g(\lambda\,\theta,|u(\lambda\,\theta)|^2)\big)
u(\lambda\,\theta)\overline{\psi_k(\theta)}\,dS(\theta).
\end{equation}
Equations (\ref{eq:lin-nonlin}) and (\ref{eq:2rad}) imply that, for every $k$,
\begin{equation*}
  -\varphi_k''(\lambda)-\frac{N-1}{\lambda}\varphi_k^\prime(\lambda)+
  \frac{\mu_k(\A,a)}{\lambda^2}\varphi_k(\lambda)=
  \zeta_k(\lambda),\quad\text{in }(0,\tilde r).
\end{equation*}
A direct calculation shows that, for some
$c_1^k(R),c_2^k(R)\in\R$,
\begin{equation}\label{eq:42}
\varphi_k(\lambda)=\lambda^{\sigma^+_k}
\bigg(c_1^k(R)+\int_\lambda^R\frac{s^{-\sigma^+_k+1}}{\sigma^+_k-\sigma^-_k}
\zeta_k(s)\,ds\bigg)+\lambda^{\sigma^-_k}
\bigg(c_2^k(R)+\int_\lambda^R\frac{s^{-\sigma^-_k+1}}{\sigma^-_k-\sigma^+_k}
\zeta_k(s)\,ds\bigg),
\end{equation}
where
\begin{equation}\label{eq:sigmapm}
  \sigma^+_k=-\frac{N-2}{2}+\sqrt{\bigg(\frac{N-2}
    {2}\bigg)^{\!\!2}+\mu_k(\A,a)}\quad\text{and}\quad
  \sigma^-_k=-\frac{N-2}{2}-\sqrt{\bigg(\frac{N-2}{2}\bigg)^{\!\!2}+\mu_k(\A,a)}.
\end{equation}
In view of Lemma \ref{l:blowup}, there exist $j_0,m\in\N$,
$j_0,m\geq 1$ such that $m$ is the multiplicity of the eigenvalue
$\mu_{j_0}(\A,a)=\mu_{j_0+1}(\A,a)=\cdots=\mu_{j_0+m-1}(\A,a)$ and
\begin{equation}\label{eq:15}
  \gamma=\lim_{r\rightarrow 0^+} {\mathcal N}(r)=\sigma_{i}^+,
  \quad i=j_0,\dots,j_0+m-1.
\end{equation}
The Parseval identity yields
\begin{equation}\label{eq:17}
H(\lambda)=\int_{{\mathbb
    S}^{N-1}}|u(\lambda\,\theta)|^2\,dS(\theta)=
\sum_{k=1}^{\infty}|\varphi_k(\lambda)|^2,\quad\text{for all
}0<\lambda\leq R.
\end{equation}
Let us assume by contradiction that $\lim_{\lambda\to
0^+}\lambda^{-2\gamma}H(\lambda)=0$.
Then, (\ref{eq:15}) and (\ref{eq:17})
imply that
\begin{equation}\label{eq:11}
\lim_{\lambda\to0^+}\lambda^{-\sigma_{i}^+}\varphi_{i}(\lambda)=0\qquad
\text{for any } i\in\{j_0,\dots,j_0+m-1\} \ .
\end{equation}
We claim that the functions
\begin{equation} \label{claim-int}
s\mapsto \frac{s^{-\sigma^+_i+1}}{\sigma^+_i-\sigma^-_i}
\zeta_i(s), \quad s\mapsto
\frac{s^{-\sigma^-_i+1}}{\sigma^-_i-\sigma^+_i} \zeta_i(s),
\end{equation}
belong to $L^1(0,R)$ for any $i\in \{j_0,\dots,j_0+m-1\}$.
To this purpose, we define
\begin{equation*}
Z_i(s)=\int_{B_s} |h(x)+g(x,|u(x)|^2)
||u(x)||\psi_i(x/|x|)|\, dx
\end{equation*}
for any $s\in (0,\tilde r)$ and for any $i\in
\{j_0,\dots,j_0+m-1\}$. We observe that $Z_i$ is an absolutely
continuous function whose derivative, defined for almost every $s\in
(0,\tilde r)$, is given by
\begin{equation*}
Z'_i(s)=s^{N-1} \int_{\SN} |h(s\theta)+g(s\,\theta,|u(s,\theta)|^2)||u(s\theta)||
\psi_i(\theta)|\, dS(\theta) \qquad \text{for a.e. } s\in
(0,\tilde r) .
\end{equation*}
Integrating by parts, we obtain
\begin{align} \label{int-parts}
\int_\lambda^R & \frac{s^{-\sigma^+_i+1}}{\sigma^+_i-\sigma^-_i}
|\zeta_i(s)|\,ds \leq
\int_\lambda^R\frac{s^{-\sigma^+_i+2-N}}{\sigma^+_i-\sigma^-_i}
Z'_i(s)\, ds \\
\notag & =\left[\frac{s^{-\sigma^+_i+2-N}}{\sigma^+_i-\sigma^-_i}
Z_i(s) \right]_\lambda^R-\int_\lambda^R
\frac{2-N-\sigma_i^+}{\sigma^+_i-\sigma^-_i}
s^{-\sigma^+_i+1-N}Z_i(s) \, ds.
\end{align}
From (\ref{eq:47}) and (\ref{F})
\begin{align}\label{Z}
  &|Z_i(s)| \leq \left(\int_{B_s} |h(x)+g(x,|u(x)|^2)| |u(x)|^2 dx
  \right)^{\!\!\frac12}
  \left(\int_{B_s} |h(x)+g(x,|u(x)|^2)| |\psi_i(x/|x|)|^2 dx \right)^{\!\!\frac12} \\
  \notag & \leq \left[\overline{C}^{-1}
    \bigg(\eta_0(s)+C_g\bigg(\frac{\omega_{N-1}}N\bigg)^{\!\!\frac2N}s^2+C_g
    \|u\|_{L^{2^*}(B_s)}^{2^*-2} \bigg) s^{N-2}
    \left(D(s)+\frac{N-2}2H(s)\right)\right]^{1/2} \times \\
  \notag \times & s^{\frac{N-2}2} \bigg[ \frac{{\eta_0(s)}}{N-2}
    \int_{\SN} \big(|\nabla_{\A}
    \psi_i(\theta)|^2-a(\theta)|\psi_i(\theta)|^2\big)
    dS(\theta)+\frac{N-2}2{\eta_0(s)} \int_{\SN} |\psi_i(\theta)|^2 dS(\theta)\\
\notag &\hskip4cm
    +\frac{C_g }{N^{2/2^*}}\|\psi_i\|_{L^{2^*}(\SN)}^2\bigg(\bigg(\frac{\omega_{N-1}}N
\bigg)^{\!\!\frac2N}s^2+ \|u\|_{L^{2^*}(B_s)}^{2^*-2}\bigg)
\bigg]^{\!\frac12} \\
  \notag & \leq\widetilde C_1(i)\sqrt{\mathcal N(s)+\frac{N-2}2 }
  \Big(\eta_0(s)+s^2+s^{\frac{2(q-2^*)}{q}}
\Big) s^{N-2}
 \sqrt{H(s)}\\
  \notag & \leq \widetilde C_1(i) \left(\sup_{(0,\tilde r/2)}
    \sqrt{\mathcal N+\frac{N-2}2 }\right) s^{N-2} \widetilde\eta(s)
  \sqrt{H(s)} \qquad \text{for all } s\in (0,\tilde r/2)
\end{align}
for some constant $\widetilde C_1(i)>0$ depending on $\overline{C}$,
$C_g$, $N$, $u$, $q$, and $\psi_i$, where
$$
\widetilde\eta(s):=\eta_0(s)+s^2+s^{\frac{2(q-2^*)}{q}}.
$$
We notice that, by assumption (\ref{hp:eta_0}),
$$
\frac{\widetilde\eta(s)}s\in L^1(0,\tilde r)
$$
and,
by Lemma \ref{gamma},
$$
\sup_{(0,\tilde r/2)} \sqrt{\mathcal N+\frac{N-2}2
}<+\infty .
$$
Inserting \eqref{Z} into \eqref{int-parts} we obtain
\begin{multline}\label{30}
\int_\lambda^R  \frac{s^{-\sigma^+_i+1}}{\sigma^+_i-\sigma^-_i}
|\zeta_i(s)|\,ds\\
 \leq \widetilde C_2(i)
\frac{\sqrt{H(R)}}{R^{\sigma_i^+}} \widetilde\eta(R)+ \widetilde C_2(i)
\frac{\sqrt{H(\lambda)}}{\lambda^{\sigma_i^+}} \widetilde\eta(\lambda)
+\widetilde C_3(i)\int_\lambda^R
\frac{\sqrt{H(s)}}{s^{\sigma_i^+}} \frac{\widetilde\eta(s)}{s} \, ds
\end{multline}
and using \eqref{hp:eta_0}, \eqref{1stest}, the integrability of
the first function in \eqref{claim-int} follows. The integrability
of the second function also follows since $\sigma_i^-<\sigma_i^+$.
Hence
\[
\lambda^{\sigma^+_{i}}
\bigg(c_1^{i}(R)+\int_\lambda^R\frac{s^{-\sigma^+_{i}+1}}{\sigma^+_{i}-\sigma^-_{i}}
\zeta_{i}(s)\,ds\bigg)=o(\lambda^{\sigma^-_{i}})\quad\text{as
}\lambda\to0^+,
\]
and then, since $\frac{u}{|x|}\in L^2(B_R,\C)$ and
$\frac{|x|^{\sigma^-_{i}}}{|x|}\not \in L^2(B_R,\C)$, we conclude
that there must be
\begin{equation} \label{c-2}
c_2^{i}(R)=-\int_0^R\frac{s^{-\sigma^-_{i}+1}}{\sigma^-_{i}-\sigma^+_{i}}
\,\zeta_{i}(s)\,ds.
\end{equation}
Using \eqref{Z} and \eqref{1stest}, we then deduce that
\begin{align}\label{eq:12}
 & \left|\lambda^{\sigma^-_{i}}
  \bigg(c_2^{i}(R)+\int_\lambda^R\frac{s^{-\sigma^-_{i}+1}}{\sigma^-_{i}-\sigma^+_{i}}
  \zeta_{i}(s)\,ds\bigg)\right|
  =\left|\lambda^{\sigma^-_{i}}
  \bigg(\int_0^\lambda
  \frac{s^{-\sigma^-_{i}+1}}{\sigma^+_{i}-\sigma^-_{i}}
  \zeta_{i}(s)\,ds\bigg)\right| \\
  \notag &
\leq \lambda^{\sigma^-_{i}} \int_0^\lambda
  \frac{s^{-\sigma^-_{i}+2-N}}{\sigma^+_{i}-\sigma^-_{i}}
Z'_i(s)\, ds =\frac{\lambda^{2-N}}{\sigma^+_i-\sigma^-_i}
Z_i(\lambda)-\lambda^{\sigma^-_{i}}\int_0^\lambda
\frac{2-N-\sigma_i^-}{\sigma^+_i-\sigma^-_i}
s^{-\sigma^-_i+1-N}Z_i(s) \, ds\\
\notag & =O\bigg(\lambda^{\sigma^+_i}\Big[ \widetilde\eta(\lambda)+
\int_0^\lambda \frac{\widetilde\eta(s)}{s} \, ds\Big]\bigg)
  =o(\lambda^{\sigma^+_{i}})
\end{align}
as $\lambda\to0^+$.  From (\ref{eq:42}), (\ref{eq:11}), and
(\ref{eq:12}), we obtain that
\begin{equation} \label{identita-c-1}
c_1^{i}(R)+\int_0^R\frac{s^{-\sigma^+_{i}+1}}{\sigma^+_{i}-\sigma^-_{i}}
\zeta_{i}(s)\,ds=0 \quad \text{for all } R\in (0,\tilde r/2) .
\end{equation}
Since $H\in C^1(0,\tilde r)$ and since we are assuming by
contradiction that $\lim_{\lambda\to 0^+} \lambda^{-2\gamma}
H(\lambda)=0$, we may select a sequence $\{R_n\}_{n\in
\N}\subset (0,\tilde r/2)$ decreasing to zero such that
\begin{equation*}
\frac{\sqrt{H(R_n)}}{R_n^\gamma}=\max_{s\in [0,R_n]}
\frac{\sqrt{H(s)}}{s^\gamma} .
\end{equation*}
Applying Lemma \ref{phi>0} with $\lambda_n=R_n$, we find
$i_0\in\{j_0,\dots,j_0+m-1\}$ such that, up to a subsequence,
\begin{equation} \label{phi>0-bis}
\lim_{n\to +\infty} \frac{\varphi_{i_0}(R_n)}{\sqrt{H(R_n)}}\neq 0.
\end{equation}
We  are now  going to reach a contradiction with
\eqref{identita-c-1} by choosing $i=i_0$, $R=R_n$ and $n\in \N$
sufficiently large. By \eqref{identita-c-1}, \eqref{30},
\eqref{phi>0-bis} and \eqref{1stest}, we have
\begin{align} \label{30-bis}
|c_1^{i_0}&(R_n)|=\left|\int_0^{R_n}
\frac{s^{-\sigma^+_{i_0}+1}}{\sigma^+_{i_0}-\sigma^-_{i_0}}
\zeta_{i_0}(s)\,ds \right|\\
&\notag \leq \widetilde C_2(i_0)
\frac{\sqrt{H(R_n)}}{R_n^{\gamma}} \widetilde\eta(R_n)+\widetilde
C_3(i_0)\int_0^{R_n}
\frac{\sqrt{H(s)}}{s^{\gamma}} \frac{\widetilde\eta(s)}{s} \, ds \\
\notag & \leq \widetilde C_2(i_0)
\left|\frac{\sqrt{H(R_n)}}{\varphi_{i_0}(R_n)}\right| \,
\left|\frac{\varphi_{i_0}(R_n)}{R_n^{\gamma}}\right|\widetilde\eta(R_n)+\widetilde
C_3(i_0) \left|\frac{\sqrt{H(R_n)}}{\varphi_{i_0}(R_n)}\right| \,
\left|\frac{\varphi_{i_0}(R_n)}{R_n^{\gamma}}\right| \int_0^{R_n}
\frac{\widetilde\eta(s)}{s} \, ds\\
\notag & =o\left(\frac{\varphi_{i_0}(R_n)}{R_n^{\gamma}}\right)
\end{align}
as $n\to+\infty$.  By \eqref{eq:42} with $k=i_0$, $R=R_n$ and
$\lambda=R_n$, we obtain
\begin{equation} \label{id-1-2}
\frac{\varphi_{i_0}(R_n)}{R_n^{\sigma_{i_0}^+}}=c_1^{i_0}(R_n)+c_2^{i_0}(R_n)
R_n^{\sigma_{i_0}^--\sigma_{i_0}^+}.
\end{equation}
By \eqref{c-2}, \eqref{Z} and \eqref{phi>0-bis} we have that
\begin{align} \label{30-sigma-2}
&|c_2^{i_0}(R_n) R_n^{\sigma_{i_0}^--\sigma_{i_0}^+}|
=R_n^{\sigma_{i_0}^--\sigma_{i_0}^+} \left| \int_0^{R_n}
\frac{s^{-\sigma_{i_0}^-+1}}{\sigma_{i_0}^--\sigma_{i_0}^+}
\zeta_{i_0}(s)\, ds\right| \\
\notag &\leq \widetilde C_2(i_0)
\frac{\sqrt{H(R_n)}}{R_n^\gamma}\widetilde\eta(R_n)+\widetilde C_4(i_0)
R_n^{\sigma_{i_0}^--\sigma_{i_0}^+} \int_0^{R_n}
\frac{\sqrt{H(s)}}{s^{\sigma_{i_0}^-}}\, \frac{\widetilde\eta(s)}s \, ds \\
\notag & =\widetilde C_2(i_0)
\left|\frac{\sqrt{H(R_n)}}{\varphi_{i_0}(R_n)}\right| \,
\left|\frac{\varphi_{i_0}(R_n)}{R_n^{\gamma}}\right|\widetilde\eta(R_n)
+\widetilde C_4(i_0) R_n^{\sigma_{i_0}^--\sigma_{i_0}^+}
\int_0^{R_n} \frac{\sqrt{H(s)}}{s^{\sigma_{i_0}^+}}
s^{\sigma_{i_0}^+-\sigma_{i_0}^-} \, \frac{\widetilde\eta(s)}s \, ds \\
\notag & \leq \widetilde C_2(i_0)
\left|\frac{\sqrt{H(R_n)}}{\varphi_{i_0}(R_n)}\right| \,
\left|\frac{\varphi_{i_0}(R_n)}{R_n^{\gamma}}\right|\widetilde\eta(R_n)
+\widetilde C_4(i_0)
\left|\frac{\sqrt{H(R_n)}}{\varphi_{i_0}(R_n)}\right| \,
\left|\frac{\varphi_{i_0}(R_n)}{R_n^{\gamma}}\right| \int_0^{R_n}
\frac{\widetilde\eta(s)}{s} \, ds\\
\notag & =o\left(\frac{\varphi_{i_0}(R_n)}{R_n^{\gamma}}\right).
\end{align}
Inserting \eqref{30-sigma-2} into \eqref{id-1-2} we obtain
\begin{equation*}
c_1^{i_0}(R_n)=\frac{\varphi_{i_0}(R_n)}{R_n^{\gamma}}+
o\left(\frac{\varphi_{i_0}(R_n)}{R_n^{\gamma}}\right)
\end{equation*}
as $n\to+\infty$, thus contradicting \eqref{30-bis}.
\end{pf}

The proof of Theorem \ref{t:asym} can be now obtained by proceeding similarly to \cite[Theorem 1.3]{FFT}
with small changes but for completeness we report it below.

\medskip

\begin{pfn}{Theorem \ref{t:asym}}
Identity (\ref{eq:20}) follows from part (i) of Lemma \ref{l:blowup}, thus
  there exists $k_0\in \N$, $k_0\geq 1$, such that
$$
\gamma:=\lim_{r\to
    0^+}{\mathcal N}_{u,h,f}(r)=-\frac{N-2}2+\sqrt{\Big(\frac{N-2}{2}
    \Big)^{\!2}+\mu_{k_0}(\A,a)}.
$$
  Let $m$ be the
  multiplicity of $\mu_{k_0}(\A,a)$, so that, for some $j_0\in\N$,
  $j_0\geq 1$, $j_0\leq k_0\leq j_0+m-1$,
  $\mu_{j_0}(\A,a)=\mu_{j_0+1}(\A,a)=\cdots=\mu_{j_0+m-1}(\A,a)$ and
  let $\{\psi_i:\,j_0\leq i\leq j_0+m-1\}$ be an $L^2({\mathbb
    S}^{N-1},\C)$-orthonormal basis for the eigenspace of $L_{\A,a}$
  associated to $\mu_{k_0}(\A,a)$.  Let $\lambda_n>0$, $n\in\N$ such that
$\lim_{n\to+\infty}\lambda_n=0$. Then, from part (ii) of Lemma
\ref{l:blowup} and Lemmas \ref{l:limite} and \ref{l:limitepositivo},
there exist a subsequence $\{\lambda_{n_k}\}_{k\in\N}$ and $m$ real
numbers $\beta_{j_0},\dots,\beta_{j_0+m-1}\in\R$ such that
$(\beta_{j_0},\beta_{j_0+1},\dots,\beta_{j_0+m-1})\neq(0,0,\dots,0)$
and
\begin{equation}\label{eq:23}
\lambda_{n_k}^{-\gamma}u(\lambda_{n_k}\theta)\to
\sum_{i=j_0}^{j_0+m-1} \beta_i\psi_{i}(\theta)\quad \text{in }
C^{1,\tau}({\mathbb S}^{N-1},\C) \quad \text{as }k\to+\infty
\end{equation}
and
\begin{equation} \label{eq:23grad} \lambda_{n_k}^{1-\gamma}\nabla
  u(\lambda_{n_k}\theta)\to \sum_{i=j_0}^{j_0+m-1}
  \beta_i(\gamma\psi_{i}(\theta)\theta+\nabla_{\SN} \psi_i(\theta))
  \quad \text{in } C^{0,\tau}({\mathbb S}^{N-1},\C^N) \quad \text{as
  }k\to+\infty
\end{equation}
for any $\tau\in(0,1)$. We now show that the $\beta_i$'s depend
neither on the sequence $\{\lambda_n\}_{n\in\N}$ nor on its
subsequence $\{\lambda_{n_k}\}_{k\in\N}$.

Let $R>0$ be such that $\overline{B}_{R}\subset\Omega$ and let
$\varphi_i$ and $\zeta_i$ as in (\ref{eq:22}). Then by
(\ref{expansion}) and (\ref{eq:23}) it follows that, for any
$i=j_0,\dots, j_0+m-1$,
\begin{equation}\label{eq:25}
\lambda_{n_k}^{-\gamma}\varphi_i(\lambda_{n_k}) =\int_{{\mathbb
    S}^{N-1}}\frac{u(\lambda_{n_k}\theta)}{\lambda_{n_k}^{\gamma}}
\overline{\psi_i(\theta)}\,dS(\theta)
\to\sum_{j=j_0}^{j_0+m-1} \beta_j\int_{{\mathbb
    S}^{N-1}}\psi_{j}(\theta)\overline{\psi_i(\theta)}\,dS(\theta)=\beta_i
\end{equation}
as $k\to+\infty$.  As showed in the proof of Lemma
\ref{l:limitepositivo}, for any $i=j_0,\dots, j_0+m-1$ and
$\lambda\in(0,R]$ we have
\begin{align}\label{eq:24}
\varphi_i(\lambda)&=\lambda^{\sigma^+_i}
\bigg(c_1^i(R)+\int_\lambda^R\frac{s^{-\sigma^+_i+1}}{\sigma^+_i-\sigma^-_i}
\zeta_i(s)\,ds\bigg)+\lambda^{\sigma^-_{i}}
  \bigg(\int_0^\lambda
  \frac{s^{-\sigma^-_i+1}}{\sigma^+_i-\sigma^-_i}
  \zeta_i(s)\,ds\bigg)\\
&\notag=\lambda^{\sigma^+_i}
\bigg(c_1^i(R)+\int_\lambda^R\frac{s^{-\sigma^+_i+1}}{\sigma^+_i-\sigma^-_i}
\zeta_i(s)\,ds\bigg)+o(\lambda^{\sigma^+_i})\quad\text{as }\lambda\to0^+,
\end{align}
for some $c_1^i(R)\in\R$, where $\sigma_i^\pm$ are as in
(\ref{eq:sigmapm}) and $\sigma_i^+=\gamma$. Choosing $\lambda=R$ in
the first line of (\ref{eq:24}), we obtain
\[
c_1^i(R)=R^{-\sigma^+_i}\varphi_i(R)-R^{\sigma^-_i-\sigma^+_i}\int_0^R
  \frac{s^{-\sigma^-_i+1}}{\sigma^+_i-\sigma^-_i}
  \zeta_i(s)\,ds.
\]
Using the last identity and letting $\lambda\to 0^+$ in
(\ref{eq:24}) it follows that
\[
\lambda^{-\gamma}\varphi_i(\lambda)\to
R^{-\sigma^+_i}\varphi_i(R)-R^{\sigma^-_i-\sigma^+_i}\int_0^R
  \frac{s^{-\sigma^-_i+1}}{\sigma^+_i-\sigma^-_i}
  \zeta_i(s)\,ds+\int_0^R\frac{s^{-\sigma^+_i+1}}{\sigma^+_i-\sigma^-_i}
\zeta_i(s)\,ds\quad\text{as }\lambda\to0^+,
\] 
and hence by (\ref{eq:25})
\begin{align*}
  \beta_i&= R^{-\gamma}\int_{{\mathbb S}^{N-1}}u(R\theta)
  \overline{\psi_{i}(\theta)}\,dS(\theta)
  \\
  &\quad-R^{-2\gamma-N+2}\int_{0}^R\frac{s^{\gamma+N-1}}{2\gamma+N-2}\bigg(
  \int_{{\mathbb S}^{N-1}}
  \big(h(s\,\theta)+g(s\theta,|u(s\theta)|^2)\big)u(s\,\theta)\overline{\psi_{i}(\theta)}\,dS(\theta)\bigg) ds\\
  &\quad +\int_{0}^R\frac{s^{1-\gamma}}{2\gamma+N-2}\bigg( \int_{{\mathbb
      S}^{N-1}}
  \big(h(s\,\theta)+g(s\theta,|u(s\theta)|^2)\big)
u(s\,\theta)\overline{\psi_{i}(\theta)}\,dS(\theta)\bigg) ds \ .
\end{align*}
We just proved that the $\beta_i$'s depend neither on the sequence
$\{\lambda_n\}_{n\in\N}$ nor on its subsequence
$\{\lambda_{n_k}\}_{k\in\N}$. This proves that the convergences in
(\ref{eq:23}) and \eqref{eq:23grad} actually hold as $\lambda\to
0^+$ thus completing the proof of the theorem.~\end{pfn}

\end{document}